\documentclass{article}
\usepackage{amsmath}
\usepackage{amsthm}
\usepackage{amssymb}

\newtheorem{theorem}{Theorem}

\newtheorem{proposition}[theorem]{Proposition}
\newtheorem{definition}[theorem]{Definition}
\newtheorem{example}[theorem]{Example}
\newtheorem{remark}[theorem]{Remark}
\newtheorem{corollary}[theorem]{Corollary}

\newcommand{\defin}[1]
   {\begin{definition} {\rm #1} \end{definition}}

\def\QED{\quad\blackslug\lower 8.5pt\null}

\newcommand{\crazy}[2]{\displaystyle{\mathop{#1}_{#2}}
\vphantom{\displaystyle{#1}}}

\begin{document}

\begin{center}

{\Large\bf 4-WEBS IN THE PLANE AND}

\vspace*{3mm}

{\Large\bf THEIR LINEARIZABILITY}

\vspace*{3mm}

  {\large Vladislav V. Goldberg}

\end{center}

\begin{abstract}
We investigate the linearizability problem for different classes
of 4-webs in the plane. In particular, we apply a recently found
in \cite{agl} the linearizability conditions for $4$-webs in the
plane to confirm that a 4-web $MW$ (Mayrhofer's web) with equal
curvature forms of its 3-subwebs and a nonconstant basic invariant
is always linearizable (this result was first obtained in
\cite{m28}); it also follows from the papers \cite{na96} and
\cite{na98}). Using the same conditions, we also prove that such a
4-web with a constant basic invariant (Nakai's web) is
linearizable if and only if it is parallelizable. We also study
four classes of the so-called almost parallelizable 4-webs $APW_a,
a = 1, 2, 3, 4$ (for them the curvature $K = 0$ and the basic
invariant is constant on the leaves of the web foliation $X_a$),
and prove that a 4-web $APW_a$ is linearizable if and only if it
coincides with a 4-web $MW_a$ of the corresponding special class
of 4-webs $MW$. The existence theorems are proved for all the
classes of 4-webs considered in the paper.
\end{abstract}

\setcounter{section}{-1}

\section{Introduction}

Let $W$ be a 4-web given by four codimension one foliations $X_a,
a = 1, 2, 3, 4,$ of curves on a two-dimensional manifold $M^2$.
The web $W$ is linearizable (rectifiable) if it is equivalent to a
linear 4-web,  i.e., a 4-web formed by  four codimension one
foliations of straight lines (see Definition 4 in subsection 2.1).

It is well-known that the geometry of a 4-web $W$ in the plane is
completely determined by the curvature $K$ of its 3-subweb [1, 2,
3], defined by its first three foliations $X_1, X_2$, and $X_3$,
and by the basic invariant $a$ and their covariant derivatives
(see \cite{g93}, or  \cite{g77}, \cite{g80},  \cite{g88}, Section
7.1, where 4-webs of codimension $r, r \geq 1,$ on $M^{2r}$ were
considered). Note that the basic invariant $a$ is the cross-ratio
of the four tangents to the leaves $F_2, F_1, F_3, F_4$ of the
foliations $X_2, X_1, X_3, X_4$ passing through the point $p$ (see
Section 3).

There are many special classes of 4-webs of codimension $r, \, r >
1,$ on $M^{2r}$ and  only a few special classes of 4-webs in the
plane (of codimension one). The reason for this is that in the
case of codimension $r = 1$, most of the classes coincide with the
class of parallelizable 4-webs, i.e., 4-webs which are equivalent
to a 4-web $W$ formed by four foliations of parallel straight
lines of an affine plane $\mathbb{A}^2$ (see Definition 3 in
subsection 2.1).

The known (different) special classes of 4-webs in the plane are:
parallelizable 4-webs, almost parallelizable 4-webs (4-webs whose
basic invariant is constant on the leaves of one of the web
foliations), Mayrhofer's 4-webs $MW$ (4-webs with equal curvature
forms of their 3-subwebs and nonconstant basic invariant) (see
\cite{g88}, Section 7.2, \cite{g93}, and \cite{g97}) and Nakai's
4-webs $NW$ (4-webs with a constant basic invariant) (see
\cite{na96} and \cite{na98}). These classes of 4-webs are
characterized by some invariant conditions imposed on the
curvature $K$, the basic invariant $a$, and their covariant
derivatives.  There is  also the class of 4-webs of maximal rank
(see \cite{h}) but the relations between $K, a$, and their
covariant derivatives for this class are unknown.

A criterion of linearizability is very important in web geometry
and in its applications. For 3-webs and 4-webs, the problem of
linearizability was posed by Blaschke (see \cite{bl}, \S\S 17 and
42). Blaschke claimed that because of complexity of calculations
involving high order differential neighborhoods, it is hopeless to
find such a criterion.

Recently the conditions of linearizability for $d$-webs, $d \geq
4$, on 2-dimensional manifolds were found by Akivis, Goldberg, and
Lychagin in \cite{agl}. For 4-webs, these conditions were
presented in two different forms: in the form of two relations
expressing the covariant derivatives $K_1$ and $K_2$ of the
curvature $K$ in terms of the curvature $K$ itself, the basic
invariant $a$ and its covariant derivatives $a_i, a_{ij}, a_{ijk}$
of the first three orders, and in the form of two partial
differential equations of the fourth order for the web functions
$z = f (x, y)$ and $u = g (x, y)$.

In the current paper we consider special classes of 4-webs in the
plane and investigate the linearizability problem for them using
the first form of linearizability conditions found in \cite{agl}.

There are two different classes of $4$-webs  with equal curvature
forms of its $3$-subwebs: 4-webs $MW$ whose basic invariant is
nonconstant and 4-webs $NW$ whose basic invariant is constant. We
prove the result which follows from the main result of Mayrhofer
\cite{m28} as well as from Nakai's results in \cite{na96} and
\cite{na98}: for 4-webs $MW$ the curvatures and the curvature
forms of all their 3-subwebs vanish, i.e., all 3-subwebs are
hexagonal (see also \cite{g93} and \cite{g97}). 4-webs all
3-subwebs of which are hexagonal were introduced by Mayrhofer (see
\cite{m28}). We arrived to Mayrhofer's 4-webs $MW$ solving the
problem of finding $4$-webs with a nonconstant basic invariant and
equal curvature forms of its $3$-subwebs. Mayrhofer was first who
proved that such 4-webs are equivalent to 4-webs all foliations of
which are pencils of straight lines. There are different proofs of
Mayrhofer's theorem: see \cite{m28}, \cite{m29}, \cite{r28},
\cite{bb}, \S 10. All these proofs are rather complicated. We
think that the reason for this is that the authors of these papers
did not use the invariant characterization of 4-webs $MW$ similar
to the conditions (42)--(45). Using the linearizability condition
of 4-webs indicated above, we present a straightforward short
proof of Mayrhofer's theorem (Theorem 10): \textit{a $4$-web $MW$
with a nonconstant basic invariant and equal curvature forms of
its $3$-subwebs is linearizable}.

Note that Mayrhofer's result on 4-webs $MW$ is stronger. He proved
that the 4-webs $MW$ are equivalent to 4-webs formed by four
pencils of straight lines. However, since in this paper we deal
only with 4-web's linearizability, we do not consider these four
pencils.

As to the 4-webs with a constant basic invariant and equal
curvature forms of its $3$-subwebs and (they were studied by Nakai
in \cite{na96} and \cite{na98}), using the same linearizability
condition, we prove that \emph{such $4$-webs are linearizable if
and only if they are parallelizable}.

We also consider four classes of the so-called almost
parallelizable 4-webs $APW_a, \, a = 1, 2, 3, 4$ (for them the
curvature $K = 0$ and the basic invariant is constant on the
leaves of one of the web foliations $X_a$). Using again the same
linearizability condition, we prove that \emph{a $4$-web $APW_a$
is linearizable if and only if it coincides with a $4$-web $MW_a$
of the corresponding special class of $4$-webs $MW$.}

The existence theorems are proved for all classes of 4-webs
considered in the paper. Note that some of the classes of 4-webs,
for example, almost parallelizable 4-webs $APW_a$ were not known
earlier.

\section{Basic Notions and Equations}

\textbf{1. 4-webs in the plane.} Let $M^2$ be a two-dimensional
analytic manifold.

\defin{\label{def:1}
A \textit{$4$-web} $W$ is given in an open domain $D$ of $M^2$ by
four foliations $X_a,\; a = 1, 2, 3, 4,$ of curves if the tangent
lines to these  curves through any point $p \in D$ are in general
position, i.e., the tangents to the curves of different foliations
are distinct and two curves of different foliations have in $D$ at
most one common point. }

All functions and differential forms which we will introduce on
$M^2$ will be real of class $C^\infty$, or holomorphic in the
complex case.

\defin{\label{def:2} Two 4-webs $W$ and $\widetilde{W}$ with
domains $D \subset M^2$ and $\widetilde{D} \subset
\widetilde{M}^2$ are \textit{equivalent} if there exists a
diffeomorphism $\phi: D \rightarrow \widetilde{D}$ (respectively
an analytic isomorphism) of the web domains such that $\phi(X_a) =
\widetilde{X_a}$.}

\textbf{2. Basic equations of a 3-web.} Suppose that  $(x,y)$ are
local coordinates in $D \subset M^2$ and  the foliations $X_a, a =
1, 2, 3, 4$,   of a 4-web  $W$ are given by the equations
\begin{equation}\label{eq:1}
  u_a(x,y) = \mbox{const}, \;\;\; a = 1, 2, 3, 4,
\end{equation}
where $u_a(x,y)$ are differentiable functions.

In order to find the basic equations of a 4-web $W$, we first
consider its 3-subweb [1, 2, 3] defined  by the foliations
$X_\alpha,  \alpha = 1, 2, 3$. If we multiply each of the total
differentials $du_\alpha$ of the functions $u_\alpha(x,y)$ from
(1) by a factor $g_\alpha(x,y) \neq 0$, we obtain the Pfaffian
forms
\begin{equation}\label{eq:2}
      \omega_\alpha = g_\alpha du_\alpha,\;\;\; \alpha = 1, 2, 3.
\end{equation}
The differential equation $\omega_\alpha = 0$ defines the
foliation $X_\alpha$ of curves. The factors $g_\alpha$ can be
chosen in such a way that for all the points $(x,y) \in D$ and for
all the directions $dx/dy$, we have the equation
\begin{equation}\label{eq:3}
              \omega_1 + \omega_2 + \omega_3 = 0
\end{equation}
(see \cite{bl}, \S 6).

It follows from (2) that
\begin{equation}\label{eq:4}
d \omega_\alpha = \omega_\alpha \wedge \theta_\alpha,
\end{equation}
where
$$
          \theta_\alpha = - d(\ln g_\alpha).
$$
Taking the exterior derivatives of (3) and using (3) and (4), we
obtain
$$
             \omega_1 \wedge (\theta_1 - \theta_3) +
             \omega_2 \wedge (\theta_2 - \theta_3) = 0.
$$
Application of Cartan's lemma to this equation implies
\begin{equation}\label{eq:5}
\left\{ \begin{array}{ll}
    \theta_1 - \theta_3 = a\omega_1 + b\omega_2, \\
    \theta_2 - \theta_3 = b\omega_1 + c\omega_2.
\end{array} \right.
\end{equation}

It follows from (5) that
$$
\theta_1 - (a - b)\omega_1 = \theta_3 + b(\omega_1 + \omega_2) =
\theta_2 - (c - b)\omega_2 = \theta
$$
or
\begin{equation}\label{eq:6}
\left\{ \begin{array}{ll}
\theta_1 = \theta + (a - b)\omega_1, \\
\theta_2 = \theta + (c - b)\omega_2, \\
\theta_3 = \theta + b\omega_3.
\end{array} \right.
\end{equation}

Equations (6) allow us to write equations (4) in the form
\begin{equation}\label{eq:7}
d\omega_\alpha =  \omega_\alpha \wedge \theta, \;\;\; \alpha = 1,
2, 3.
\end{equation}
The form $\theta$ is called the \emph{connection form} of the
three-web [1, 2, 3].

Exterior differentiation of (7) leads to the cubic exterior
equation $d\theta \wedge \omega_\alpha = 0$
  from which we easily conclude that
\begin{equation}\label{eq:8}
         d\theta = K\Omega,
\end{equation}
where $K$ is the {\em curvature} and
\begin{equation}\label{eq:9}
\Omega = \omega_1 \wedge \omega_2
        = \omega_2 \wedge \omega_3
        = \omega_3 \wedge \omega_1
\end{equation}
is the {\em surface element} of the three-web [1, 2, 3]. The form
$\Theta = K \Omega$ is called the {\em curvature form} of the
three-web [1, 2, 3]. The equations (8) and (9) are the
\emph{structure equations} of the 3-subweb [1, 2, 3] (see
\cite{bl}, \S 8). Equations (7) and (8) prove that the
\emph{affine connection} $\gamma_{123}$ is associated with
3-subweb [1, 2, 3], and $\theta$ and $K$ are the \emph{connection
form} and the \emph{curvature} of this connection.

\textbf{3. Basic equations of a 4-web.} The fourth foliation of
curves of a 4-web $W$ can be given by the equation $\omega_4 = 0$,
where
\begin{equation}\label{eq:10}
- \omega_4 = a \omega_1 + \omega_2 = 0
\end{equation}
(see \cite{g88}, Section 7.1, or \cite{g93}).

The quantity $a$ in equation (10) is called the \emph{basic
invariant} of the 4-web $W$. It satisfies the conditions
\begin{equation}\label{eq:11}
  a \neq 0, 1
\end{equation}
(see   \cite{g93}, or \cite{g77}, \cite{g80}, \cite{g88}, Section
7.1).

Equations $\omega_{a} = 0, \, a = 1, 2, 3, 4,$ defining the
foliations $X_a$, are preserved under the following concordant
transformations of the form $\omega_a$:
$$
{}^\prime \omega_a = s \, \omega_a,
$$
and only under such transformations. Under such transformations
the basic invariant $a$ is not changed, ${}'a = a$, i.e., $a$ is
an absolute invariant.

It is easy to show (see  \cite{g93}, or \cite{g77}, \cite{g80},
\cite{g88}, Section 7.1) that \emph{at a point $p \in M^2$, the
basic invariant $a$ of the web $W$ equals the cross-ratio of the
four tangents to the leaves $F_2, F_1, F_3, F_4$ of the foliations
$X_2, X_1, X_3, X_4$ passing through the point $p$.}

The basic invariant $a$ of the 4-web $W$ is an absolute invariant
of $W$. It satisfies the differential equation
\begin{equation}\label{eq:12}
       da = a_1 \, \omega_1 + a_2 \,\omega_2,
\end{equation}
where $a_i \in C^{\infty}(M^2), \;i = 1, 2$ are the first
covariant derivatives of $a$.

\textbf{4. Prolongations of the basic equations.} In what follows,
we will need  prolongations of equations (8). Exterior
differentiation of (8) by means of (7) gives the following
exterior cubic equation:
\begin{equation}\label{eq:13}
          \nabla K \wedge \omega_1 \wedge \omega_2 = 0,
\end{equation}
where $\nabla K = dK - 2K\theta$. It follows from (13) that
\begin{equation}\label{eq:14}
      \nabla K = K_1 \omega_1 + K_2 \omega_2,
\end{equation}
where $K_i \in C^{\infty}(M^2), \;\;i = 1,2$. Equation (14) shows
that the curvature $K$ is a relative invariant of weight two.

We recall that a quantity $u$ is a relative invariant of weight
$k$ if
$$
du - k \theta = u_1 \omega_1 + u_2 \omega_2
$$
(see \cite{ag}, p. 36) If $k = 0$, then $u$ is an absolute
invariant.

We will also need two prolongations of equations (12). Taking
exterior derivatives of equations (12), we obtain the following
exterior quadratic equation:
\begin{equation}\label{eq:15}
       \nabla a_1 \wedge \omega_1 + \nabla a_2 \wedge \omega_2 = 0,
\end{equation}
where $\nabla a_i = d a_i - a_i \theta, \, i = 1, 2$. Applying
Cartan's lemma to equation (15), we find that
\begin{equation}\label{eq:16}
       \nabla a_1 = a_{11} \omega_1 + a_{12} \omega_2, \;\;
\nabla a_2 = a_{12} \omega_1 + a_{22} \omega_2,
\end{equation}
where  $a_{ij} \in C^{\infty}(M^2), \; i, j = 1,2$.

Taking exterior derivatives of equations (16), we obtain the
following exterior quadratic equations:
\begin{equation}\label{eq:17}
\renewcommand{\arraystretch}{1.5}
\left\{
\begin{array}{ll}
\nabla a_{11} \wedge \omega_1 + \nabla a_{12} \wedge \omega_2 +
a_1 K \omega_1 \wedge \omega_2= 0,\\ \nabla a_{12} \wedge \omega_1
+ \nabla a_{22} \wedge \omega_2 +
  a_2 K \omega_1 \wedge \omega_2 = 0,
\end{array}
\right.
\renewcommand{\arraystretch}{1}
\end{equation}
where $\nabla a_{ij} = da_{ij} - 2 a_{ij} \theta, \, i, j = 1, 2$.
Applying Cartan's lemma to equations (17), we find that
\begin{equation}\label{eq:18}
\renewcommand{\arraystretch}{2.3}
\left\{
\begin{array}{ll}
\nabla a_{11} = a_{111} \omega_1 + (a_{112} + \displaystyle \frac{2}{3} a_1 
K) \omega_2, \\
\nabla a_{12} = (a_{112} -  \displaystyle \frac{1}{3} a_1 K)
\omega_1
    + (a_{122} +  \displaystyle \frac{1}{3} a_2 K) \omega_2,\\
\nabla a_{22} = (a_{122} -  \displaystyle \frac{2}{3} a_2 K)
\omega_1 + a_{222} \omega_2,
\end{array}
\right.
\renewcommand{\arraystretch}{1}
\end{equation}
where  $a_{ijk} \in C^{\infty}(M^2), \; i, j, k = 1,2$.

It follows from (16) and (18) that $a_i$ and $a_{ij}$ are relative
invariants of weight one and two, respectively.

\textbf{5. Connection forms, curvatures, and  curvature forms of
3-subwebs of a 4-web.} Next we will indicate the formulas for the
connection forms, the curvatures, and  the curvature forms of all
3-subwebs [1, 2, 3] [1, 2, 4], [1, 3, 4] and [2, 3, 4] of a 4-web
$W$.

For the 3-subweb [1, 2, 3], the indicated quantities are:
\begin{equation}\label{eq:19}
  \crazy{\theta}{123} = \theta,  \;\;
  \crazy{K}{123} = K, \;\;  \crazy{\Theta}{123} = K
\omega_1 \wedge \omega_2.
\end{equation}

As was shown in \cite{g93} (see also \cite{g88}, Section 7.1), for
the 3-subweb [1, 2, 4] these quantities are
\begin{equation}\label{eq:20}
  \renewcommand{\arraystretch}{1.8}
\begin{array}{ll}
\crazy{\theta}{124} = \theta - \frac{a_2}{a}
\omega_2, \\
  \crazy{K}{124} = \frac{1}{a}
\Bigl(K - \frac{a_{12}}{a}
+  \frac{a_{1} a_2}{a^2}\Bigr), \\
  \crazy{\Theta}{124}  = a \crazy{K}{124}
\omega_1 \wedge \omega_2;
\end{array}
  \renewcommand{\arraystretch}{1}
\end{equation}
for the 3-subweb [1, 3, 4], they  are:
\begin{equation}\label{eq:21}
  \renewcommand{\arraystretch}{1.8}
\begin{array}{ll}
\crazy{\theta}{134} = \theta + \frac{a_2}{1 - a}
(\omega_1 + \omega_2), \\
  \crazy{K}{134} = \frac{1}{a - 1}
\Bigl[K + \frac{a_{2}(a_1 - a_2)}{(1 - a)^2}
+  \frac{a_{12} - a_{22}}{1 - a}\Bigr], \\
  \crazy{\Theta}{134}  = (a - 1)  \crazy{K}{134}
\omega_1 \wedge \omega_2;
\end{array}
  \renewcommand{\arraystretch}{1}
\end{equation}
and for the 3-subweb [2, 3, 4], they are
\begin{equation}\label{eq:22}
\renewcommand{\arraystretch}{1.8}
\begin{array}{ll}
\crazy{\theta}{234} = \theta + \frac{a_1 - a_2}{a} \omega_2
+ \frac{a_1}{1 - a} (\omega_1 + \omega_2) , \\
  \crazy{K}{234} = \frac{1}{a(a - 1)}
\Bigl[K + \frac{(2a - 1)a_1 (a_1 - a_2)}{a^2 (1 - a)^2} +
\frac{a_{11}
- a_{12}}{a(1 - a)}\Bigr], \\
  \crazy{\Theta}{234}  = a (a - 1)
  \crazy{K}{234} \omega_1 \wedge \omega_2.
\end{array}
  \renewcommand{\arraystretch}{1}
\end{equation}

\section{Special Classes of 4-Webs in the Plane}

\textbf{1. Parallelizable and linearizable  4-webs.}

\defin{\label{def:3} A 4-web $W$ formed by four
foliations of parallel straight lines in an affine plane
$\mathbb{A}^2$ is said to be a \emph{parallel $4$-web}. A 4-web
which is  equivalent to the parallel 4-web is called {\em
parallelizable}. }

For a 3-web in $M^2$ the definitions of parallel and
parallelizable 3-webs are similar.

\textbf{Remark}. A parallel 4-web defined by the foliations
$X_b,\, b = 1, 2, 3, 4$, can be given by the functions
$$
u_1 = x, \;\; u_2 = y, \;\; u_3 = x + y, \;\; u_4 = a x + y
$$
  (cf.  \cite{g93} for parallel 3-webs).

\setcounter{theorem}{0}
\begin{theorem} A $4$-web $W$ is parallelizable if and only if
the curvature $K$ vanishes, and the basic invariant is constant in
$M^2$, i.e., if and only if the following conditions are
satisfied:
\begin{equation}\label{eq:23}
K = 0, \;\;\; a_1 = 0, \;\;\; a_2 = 0.
\end{equation}
\end{theorem}

\begin{proof} The condition $K = 0$ implies that the leaves of
the foliations $X_1, X_2$, and $X_3$ can be defined as level sets
of the functions $u_1, u_2$, and $u_3$ in Remark above. The
conditions $a_1 = 0, \, a_2 = 0$ imply that the basic invariant
$a$ is constant. Hence leaves of the foliation $X_4$ are parallel
straight lines which are level sets of the function $u_4$ in
Remark. The converse is obvious. \end{proof}

For another proof see \cite{g77} or \cite{g80}, or \cite{g88},
Section 7.2.

\setcounter{theorem}{3}

\defin{\label{def:4} A 4-web $W$ formed by 4 foliations
of straight lines is said to be \emph{linear}. A 4-web which is
equivalent to a linear 4-web is called \emph{linearizable}
(\emph{rectifiable}). }

It follows from Definitions 3 and 4 that a parallelizable 4-web is
linearizable but the converse is not true.

In \cite{agl}  the following criterion of linearizability for
4-webs in the plane was found. \setcounter{theorem}{1}
\begin{theorem}
A $4$-web $W$ is linearizable if and only if the curvature $K$,
its covariant derivatives $K_1$ and $K_2$, the basic invariant
$a$, and its covariant derivatives $a_i, a_{ij}, a_{ijk}$ of the
first three orders satisfy the following equations:
\begin{equation}\label{eq:24}
\renewcommand{\arraystretch}{1.5}
\begin{array}{ll}
    K_1  =  &\!\!\!\! \displaystyle\frac{1}{a-a^2}
    \Bigl[\displaystyle \frac{1}{3}\Bigl(a_1 (1-a)
    +  aa_2\Bigr)K -  a_{111} + (2+a) a_{112}
   - 2 a a_{122}\Bigr]  \\
  &\!\!\!\!+ \displaystyle  \frac{1}{(a-a^2)^2}
\Biggl\{\Bigr[(4-6a)a_1 +(-2+3a+a^2)a_2\Bigl] a_{11}\\
  &\!\!\!\!
+\Bigl[(-6+7a+2a^2)a_1 + (2a-3a^2)a_2\Bigr]a_{12}
+ \Bigl[(2a(1-a)a_1 - 2a^2 a_2\Bigr]a_{22}\Biggr\}\\
&\!\!\!\!+ \displaystyle  \frac{1}{(a-a^2)^3}
\Biggl\{(-3+8a-6a^2)(a_1)^3 -2a^3(a_2)^3 \\
&\!\!\!\!+(6 - 15a+9a^2+2a^3)(a_1)^2 a_2 +(-2a+6a^2-3a^3)a_1
(a_2)^2\Biggr\},
   \end{array}
\renewcommand{\arraystretch}{1}
\end{equation}
\begin{equation}\label{eq:25}
\renewcommand{\arraystretch}{1.5}
\begin{array}{ll}
   K_2  =  &\!\!\!\! \displaystyle\frac{1}{a-a^2}
   \Bigl[\displaystyle \frac{1}{3} \Bigl(a_1
   + a_2 (a - 1)\Bigr) K + 2 a_{112} - (2a+1) a_{122}
   + a a_{222}\Bigr]   \\
  &\!\!\!\!+ \displaystyle  \frac{1}{(a-a^2)^2}
\Biggl\{\Bigr[(2a_1 +(2a-2)a_2\Bigl] a_{11} + \Bigl[(-5 +6a)a_1
+ (2-3a - 2 a^2)a_2\Bigr]a_{12}\\
&\!\!\!\!+ \Bigl[(1-a-2a^2)a_1+ 2a^2 a_2\Bigr]a_{22}\Biggr\} +
\displaystyle  \frac{1}{(a-a^2)^3}
\Biggl\{(4a-2)(a_1)^3 + a^3(a_2)^3 \\
&\!\!\!\!+(5 - 12a+6a^2)(a_1)^2 a_2 + (-2 + 5a- 3a^2-2a^3)a_1
(a_2)^2\Biggr\}.
   \end{array}
\renewcommand{\arraystretch}{1}
\end{equation}
\end{theorem}

\begin{proposition} Linearizable $4$-webs exist, and the general solution
of the system of equations defining such webs depends on four
arbitrary  functions of one variable.
\end{proposition}

\begin{proof}
Taking exterior derivatives of (18) and applying (4), (8), (12),
(14), (16), (18), (24), and (25), we get three exterior quadratic
equations
\begin{equation}\label{eq:26}
\renewcommand{\arraystretch}{1.5}
\begin{array}{ll}
\nabla a_{111} \wedge \omega_1 + \nabla a_{112} \wedge \omega_2 +
(\ldots) \omega_1 \wedge \omega_2= 0,\\
\nabla a_{112} \wedge \omega_1 + \nabla a_{122} \wedge \omega_2 +
(\ldots) \omega_1 \wedge \omega_2= 0,\\
\nabla a_{122} \wedge \omega_1 + \nabla a_{222} \wedge \omega_2 +
(\ldots) \omega_1 \wedge \omega_2= 0,
   \end{array}
\renewcommand{\arraystretch}{1}
\end{equation}
where the coefficients $(\ldots)$ depend on $K, a, a_i$, and
$a_{ij}$. If we substitute the expressions for $K_1$ and $K_2$
from (24) and (25) into equation (14) and take exterior
derivatives of the resulting equation using the same (4), (8),
(12), (14), (16), (18), (24), and (25), we get one more exterior
quadratic equation of the form
\begin{equation}\label{eq:27}
\renewcommand{\arraystretch}{1.5}
\begin{array}{ll}
&\!\!\!\! [- \nabla a_{111} + (2 + a) \nabla a_{112}
- 2 a \nabla a_{122}]  \wedge \omega_1\\
+ &\!\!\!\! [2 \nabla a_{112} - (2a + 1)  \nabla a_{122} + a
\nabla a_{222}]  \wedge \omega_2 + (\ldots) \omega_1 \wedge
\omega_2= 0,
   \end{array}
\renewcommand{\arraystretch}{1}
\end{equation}
where the coefficients $(\ldots)$ again depend on $K, a, a_i$, and
$a_{ij}$.

The system defining linearizable 4-webs consists of the Pfaffian
equation (14) (in which $K_1$ and $K_2$ are replaced by their
values (24) and (25)), (3), (16), (18), and exterior quadratic
equations (26), (27). For this system we have 4 unknown forms
$\nabla a_{ijk}$. Hence $q = 4$. Since the  exterior quadratic
equations (26) and (27) are obviously independent, and their
number is four, we have $s_1 = 4$. Thus $s_2 = q - s_1 = 0$, and
the Cartan number $Q = s_1 + 2 s_2 = 4$. On the other hand, if we
apply Cartan's lemma to (26), we find the forms $\nabla a_{ijk}$
in terms of $\omega_i$, and there will be 5 coefficients in the
expansions of these forms. Substituting these expansions into
(27), we find one of these coefficients. So, the dimension $N$ of
the space of integral elements over a point equals four,  $N = 4$.
  Thus $Q = N$, and by Cartan's test (see \cite{bcggg}, p. 120),
  the system is in involution, and its general solution depends
  on four functions of one variable.
\end{proof}

\textbf{2. 4-webs whose basic invariant is constant on the leaves
of one of web foliations.}

\begin{theorem}  The basic invariant $a$ of a $4$-web is  constant
on the leaves of the web foliation $X_1, X_2, X_3,$ and $X_4$ if
and only if there is respectively the following relation between
the covariant derivatives $a_1$ and $a_2$ of this invariant:
\begin{equation}\label{eq:28}
  a_2 = 0,
\end{equation}
\begin{equation}\label{eq:29}
  a_1 = 0,
\end{equation}
\begin{equation}\label{eq:30}
  a_1 = a_2,
\end{equation}
\begin{equation}\label{eq:31}
  a_1 = a a_2.
\end{equation}
\end{theorem}

\begin{proof} In fact, it follows from (1.3) that
conditions (29)--(31) are necessary and sufficient for the
following congruences:
$$
\renewcommand{\arraystretch}{1.5}
\begin{array}{ll}
da \equiv 0 \pmod{\omega_1}, \\
da \equiv 0 \pmod{\omega_2},\\
da \equiv 0 \pmod{\omega_1 + \omega_2},\\
da \equiv 0 \pmod{a\omega_1 + \omega_2},
   \end{array}
\renewcommand{\arraystretch}{1}
$$
respectively.
\end{proof}

\begin{proposition}  $4$-webs whose basic invariant is constant on the
leaves of one of web foliations exist, and the general solution of
the system of equations defining such webs depends on one
arbitrary function of two variables.
\end{proposition}

\begin{proof} Suppose, for example, that the  basic invariant is
constant on the leaves of the foliation $X_2$ defined by the
equation $\omega_2 = 0$. Then condition (29) holds. For such a
4-web, we have the Pfaffian equations (12) and (16). By (29), the
former becomes
\begin{equation}\label{eq:32}
da = a_2 \omega_2.
\end{equation}
Taking the exterior derivatives of these two Pfaffian equations,
we arrive at the following two exterior quadratic equations:
\begin{equation}\label{eq:33}
\renewcommand{\arraystretch}{1.5}
\begin{array}{ll}
\nabla K_1 \wedge \omega_1
+ \nabla K_2 \wedge \omega_2 = 0, \\
\nabla a_2 \wedge \omega_2 = 0,
   \end{array}
\renewcommand{\arraystretch}{1}
\end{equation}
where $\nabla K_i = d K_i - 3 K_i \theta, \, \nabla a_2 = d a_2 -
a_2 \theta$. Thus, we have three unknown functions, and $q = 3$.
Since there are two exterior quadratic equations, we have $s_1 =
2$. As a result, $s_2 = q - s_1 = 1$, and the Cartan number $Q =
s_1 + 2 s_2 = 4$. Applying Cartan's lemma to the exterior
quadratic equations (33), we find that
$$
\renewcommand{\arraystretch}{1.5}
\begin{array}{ll}
\nabla K_1  = K_{11}  \omega_1 + K_{12} \omega_2, \\
\nabla K_2  = K_{12}  \omega_1 + K_{22} \omega_2, \\
\nabla a_2  = a_{22}  \omega_2.
   \end{array}
\renewcommand{\arraystretch}{1}
$$
Hence  the dimension $N$ of the space of integral elements over a
point equals 4,  $N = 4$. Therefore, $Q=N$, and by Cartan's test,
our system of equations is in involution, and its general solution
depends on one arbitrary function of two variables.

The proof in other cases is similar. The difference is that in the
cases when a 4-web satisfies equations (28), (30), and (31),
equation (32) has respectively the form
$$
da = a_1 \omega_1,
$$
$$
da = a_1 (\omega_1 + \omega_2),
$$
and
$$
da = a_2 (a\omega_1 + \omega_2),
$$
and the last equation (33) becomes
$$
\nabla a_1 \wedge \omega_1 = 0,
$$
$$
\nabla a_1 \wedge (\omega_1 + \omega_2) = 0,
$$
and
$$
\nabla a_2 \wedge (a\omega_1 + \omega_2) = 0,
$$
respectively.
\end{proof}

\textbf{3. Almost parallelizable 4-webs.} \setcounter{theorem}{4}
\defin{\label{def:5} We say that a 4-web $W$ is
  \emph{almost parallellizable} if its $3$-subweb $[1, 2, 3]$
is parallelizable (i.e.,  $K = 0$), and its basic invariant is
constant on the leaves of one of the web foliations. }

So, according to Theorem 4, there are four classes of almost
parallelizable 4-webs characterized by the following conditions:
\begin{equation}\label{eq:34}
  K = 0, \;\; a_2 = 0,
\end{equation}
\begin{equation}\label{eq:35}
  K = 0, \;\; a_1 = 0,
\end{equation}
\begin{equation}\label{eq:36}
  K = 0, \;\; a_1 = a_2,
\end{equation}
\begin{equation}\label{eq:37}
  K = 0, \;\; a_1 = a a_2.
\end{equation}

We denote the 4-webs of these four classes by $APW_a, \, a = 1, 2,
3, 4$.

\setcounter{theorem}{5}
\begin{proposition}  The almost parallelizable $4$-webs
$APW_a$ exist, and the general solution of the system of equations
defining such webs depends on one arbitrary function of one
variable.
\end{proposition}

\begin{proof}
In fact, in each of the cases (34)--(37), the condition $K = 0$
implies  $K_1 = 0, \, K_2 = 0$. Thus, we have $q = 1$. For a 4-web
$APW_1$, in (33) (and in similar systems corresponding to 4-webs
$APW_2, APW_3$, and $APW_4$) only the second equation remains.
Hence Cartan's characters are $s_1 = 1$ and $s_2 = 0$. It is easy
to see that in this case $Q = N = 1$, and by Cartan's test, our
system of equations is involutive, and its general solution
depends on one arbitrary function of one variable.
\end{proof}

The following theorem establishes the values of the curvatures
$\crazy{K}{abc}$ for the 4-webs $APW_a$.

\begin{theorem} The curvatures
$\crazy{K}{a b c}$ for the $4$-webs $APW_a$ are
\begin{description}
\item[(i)] For the $4$-webs $APW_1$:
\begin{equation}\label{eq:38}
\crazy{K}{123} = \crazy{K}{124} = \crazy{K}{134}= 0, \;\;\;
\crazy{K}{234} = \frac{(1-2a) a_1^2 - (a-a^2)a_{11}}{(a - a^2)^3};
\end{equation}
\item[(ii)] For the $4$-webs $APW_2$:
\begin{equation}\label{eq:39}
\crazy{K}{123} = \crazy{K}{124} = \crazy{K}{234}= 0, \;\;\;
\crazy{K}{134} = \frac{a_2^2 + (1-a)a_{22}}{(1 - a)^3};
\end{equation}
\item[(iii)] For the $4$-webs $APW_3$:
\begin{equation}\label{eq:40}
\crazy{K}{123} = \crazy{K}{134} = \crazy{K}{234}= 0, \;\;\;
\crazy{K}{124} = \frac{a_1^2 - a a_{11}}{a^3};
\end{equation}
\item[(iv)] For the $4$-webs $APW_4$:
\begin{equation}\label{eq:41}
\crazy{K}{123} = 0, \;\;\; \crazy{K}{124} = - \frac{a
a_{22}}{a^2},\;\;\; \crazy{K}{134} = \frac{a_{22}}{1 - a}, \;\;\;
\crazy{K}{234}= \frac{a_{22}}{a - a^2}.
\end{equation}
\end{description}
\end{theorem}

\begin{proof} The proof is straightforward: it follows from formulas
(20)---(22) and conditions (34)---(37).
\end{proof}

\textbf{4. 4-webs with a nonconstant basic invariant and equal
curvature forms of their 3-subwebs.} It is important to emphasize
that we assume here that \emph{the basic invariant of a $4$-web is
not constant on $M_2$}, i.e., its covariant derivatives $a_1$ and
$a_2$ do not vanish simultaneously. This assumption shows that the
class of 4-webs we are going to consider is completely different
from the class considered by Nakai in \cite{na96} and \cite{na98}.
Nakai also considered 4-webs with equal curvature forms of their
3-subwebs. Assuming that these equal curvature forms of the
3-subwebs are nonvanishing, he proved that the basic invariant $a$
of a 4-web is constant on $M_2$, i.e., $a_1 = a_2 = 0$. Thus in
\cite{na96} and \cite{na98}, Nakai considered 4-webs with a
constant basic invariant. We will consider this kind of 4-webs in
subsection 6.

\begin{theorem} A $4$-web $W$ has a nonconstant basic invariant $a$
and equal curvature forms of its $3$-subwebs $[a, b, c]$ if and
only if all $3$-subwebs $[a, b, c]$ are parallelizable, i.e., if
\begin{equation}\label{eq:42}
\crazy{K}{abc} = 0, \;\;\; \crazy{\Theta}{abc} = 0, \;\;\; a, b, c
= 1, 2, 3, 4,
\end{equation}
and the second covariant derivatives $a_{ij}$ of the basic
invariant $a$ are expressed in terms of the invariant $a$ itself
and its first covariant derivatives $a_i$ as follows:
\begin{equation}\label{eq:43}
a_{11} = \frac{a_1 [(1 - 2 a) a_1 + a a_2]}{a - a^2},
\end{equation}
\begin{equation}\label{eq:44}
a_{12} = \frac{a_1 a_2}{a},
\end{equation}
\begin{equation}\label{eq:45}
a_{22} = \frac{a_2 (a_1 - a a_2)}{a - a^2}.
\end{equation}
\end{theorem}

\begin{proof}
In fact, by equation (20), it follows that if $a \neq 0$, then the
condition $\crazy{\Theta}{124} = \crazy{\Theta}{123}$ holds if and
only if equation (44) is valid.

In a similar manner, by (21) and (22), it follows that if $a \neq
0$, then the conditions $\crazy{\Theta}{134} =
\crazy{\Theta}{123}$ and $\crazy{\Theta}{234} =
\crazy{\Theta}{123}$ hold if and only if equations (45) and (43)
are valid, respectively.

To prove (42), we substitute the values of $a_{ij}$ from
(43)--(45) into equations (16). This gives
\begin{equation}\label{eq:46}
\renewcommand{\arraystretch}{1.5}
\begin{array}{ll}
d a_{1} = a_1 \theta + \displaystyle \frac{\bigl[(1 - 2 a) a_{1} +
a a_{2}\bigr] a_1}{a - a^2}
  \omega_1 + \frac{a_1 a_{2}}{a} \omega_2, \\
d a_{2} = a_2 \theta + \displaystyle\frac{1}{a - a^2} [a_{1} a_2
(1-a) \omega_1 + a_2 (a_1 - a a_2) \omega_2].
\end{array}
\renewcommand{\arraystretch}{1}
\end{equation}
Exterior differentiation of either of equations (46) implies $K =
0$. The latter condition and equations (43)--(45) lead to (42).

Thus, each of 3-subwebs $[a, b, c]$ is parallelizable. However, by
Theorem 1, this does not imply the parallelizability of the 4-web
$W$---the latter web is parallelizable if and only if all
3-subwebs $[a, b, c]$ are parallelizable and the basic invariant
$a$ is constant in the connection $\gamma_{123}$, i.e., if and
only if we have conditions (42) and $a_1 = a_2 = 0$.
\end{proof}

We  denote the 4-webs with equal curvature forms of their
3-subwebs and a nonconstant basic invariant by $MW$.

Note that conditions (42) were first obtained by Nakai (see
\cite{na96} and \cite{na98}). As to conditions (43)--(45)
characterizing the webs $MW$, they appeared here at the first
time.

4-webs all 3-subwebs of which are parallelizable (hexagonal) were
introduced by Mayrhofer (see \cite{m28}). This is the reason that
we denoted them by $MW$ (Mayrhofer's 4-webs). We arrived to
Mayrhofer's 4-webs $MW$ solving the problem of finding $4$-webs
with a nonconstant basic invariant and equal curvature forms of
its $3$-subwebs. Mayrhofer was first who proved that such 4-webs
are equivalent to 4-webs all foliations of which are pencils of
straight lines. There are different proofs of Mayrhofer's theorem:
see \cite{m28}, \cite{m29}, \cite{r28}, \cite{bb}, \S 10. However,
all these proofs are rather complicated. We think that the reason
for this is that the authors of these papers did not use the
invariant characterization of 4-webs $MW$ similar to the
conditions (42)--(45).

Now we prove the existence theorem for webs $NW$.

\begin{proposition}  The $4$-webs $MW$  exist. The system of
Pfaffian equations defining such webs is completely integrable,
and its general solution depends on three arbitrary constants.
\end{proposition}

\begin{proof} For a web $MW$, we have Pfaffian equations (12) and
(16), where the $a_{ij}$ are expressed by formulas (43)--(45). If
we take exterior derivatives of equations (16) in which $a_{ij}$
are replaced by their values (43)--(45), we arrive at the
identities. So, for the system in question, we have $s_0 = 3,\,
s_1 = 0$. The system is completely integrable, and its general
solution depends on three arbitrary constants.
\end{proof}

In the next theorem we present a straightforward short proof of
the fact that $4$-webs $MW$ are always linearizable by applying a
linearizability condition of 4-webs found in \cite{agl}. This
theorem gives the first nontrivial application of linearizability
conditions (24) and (25) for 4-webs.

\begin{theorem} $($\rm{\textbf{Mayrhofer}} \rm{\cite{m28}}$)$
A $4$-web $MW$  is linearizable.
\end{theorem}

\begin{proof}
To prove this result, we must prove that conditions (24) and (25)
of linearizability are satisfied identically for the $4$-webs
$MW$. To check these conditions, we substitute the values  of the
invariants $a_{ij}$ and $a_{ijk}$ into them. As to $a_{ij}$, their
values are given by (43)--(45).

In order to find $a_{ijk}$, we differentiate (43)--(45) by using
(12) and (16) and compare the result with equations (18). Equating
the coefficients of $\omega_1$ and $\omega_2$ in the resulting
equations, we find that
\begin{equation}\label{eq:47}
\renewcommand{\arraystretch}{2.2}
\begin{array}{ll}
a_{111} = \displaystyle \frac{1}{(a - a^2)^2}
\Bigl[(6a^2 - 6a + 1) a_1^3 + (4a - 6 a^2) a_1^2 a_2 + a^2 a_1 a_2^2\Bigr],\\
a_{112} = \displaystyle \frac{1}{(a - a^2)^2}
\Bigl[(1-2a)(1-a) a_1^2 a_2 + (a - a^2) a_1 a_2^2\Bigr] = \frac{a_2 
a_{11}}{a},\\
a_{122} = \displaystyle \frac{a_1 a_2 (a_1 - a a_2)}{a^2 (1 - a)} = 
\frac{a_1 a_{22}}{a},\\
a_{222} = \displaystyle \frac{1}{(a - a^2)^2} \Bigl(a_1^2 a - 2 a
a_1 a_2^2 + a^2 a_2^3\Bigr).
\end{array}
\renewcommand{\arraystretch}{1}
\end{equation}

Next, by (42), for 4-webs $MW$, we have $K = 0$. This and (14)
imply  $K_1 = K_2 = 0$. Substituting the values of $a_{ij}$ from
(43)--(45) and the values of $a_{ijk}$ from (47) into (24) and
(25), after lengthy straightforward calculations, we come to the
identities: the coefficients of $a_1^3,\, a_1^2 a_2, \,a_1 a_2^2$
and $a_2^3$ vanish. Hence any 4-web $MW$ is linearizable.
\end{proof}

\textbf{5. Special classes of 4-webs $\boldsymbol{MW}$.} Although
for 4-webs $MW$ the basic invariant $a$ is not constant on $M^2$,
it can be constant on the leaves of one of the foliations $X_a$ of
a web $MW$. If this is the case, we have four special classes of
4-webs $MW$. These four classes are intersections of four classes
of almost parallelizable 4-webs $APW_a$ and the class of 4-webs
$MW$. We denote webs of these classes by $MW_a$.

\begin{theorem}
The $4$-webs $MW_a$  are characterized by the following
conditions:
\begin{equation}\label{eq:48}
\renewcommand{\arraystretch}{1.8}
\begin{array}{ll}
  {\rm  Webs} \;\; MW_1: & K = 0, \;\; a_2 = 0, \;\;
  a_{11} = \displaystyle \frac{(1 - 2 a) a_1^2}{a - a^2}, \;\; a_{12} = 0, 
\;\; a_{22} = 0;\\
{\rm Webs} \; MW_2:  &  K = 0, \;\; a_1 = 0, \;\;
  a_{11} = 0, \;\; a_{12} = 0, \;\; a_{22} = -\displaystyle \frac{a_2^2 }{1 
- a};\\
   {\rm  Webs} \;\; MW_3: &  K = 0, \;\; a_1 = a_2,\;\;
a_{11} = a_{12} = a_{22} = \displaystyle\frac{a_1^2}{a};\\
  {\rm Webs} \;\; MW_4: &  K = 0, \;\; a_1 = a a_2, \;\;
a_{11} = 2 a a_2^2, \;\; a_{12} = a_2^2, \;\; a_{22} = 0.
\renewcommand{\arraystretch}{1}
\end{array}
\end{equation}
\end{theorem}

\begin{proof}
The proof follows from (34)--(37) and (43)--(45).
\end{proof}

Note that by (48) and (38)--(41), for any 4-web $MW_a$, all the
curvatures $\crazy{K}{abc}$ vanish (as it should be according to
Theorem 8 (see (42)).

\begin{proposition}  The $4$-webs $MW_a, a = 1, 2, 3, 4,$  exist. The system of
Pfaffian equations defining such webs is completely integrable,
and its general solution depends on two arbitrary constants.
\end{proposition}
\begin{proof}
In fact, by (48), in the system defining each of the webs $MW_a$
(see proof of Theorem 9), only two Pfaffian equations are
independent.
\end{proof}

The next theorem gives a relation between the 4-webs $APW_a$ and
$MW_a$.

\begin{theorem}  A $4$-web $APW_a, a = 1, 2, 3, 4,$
is a $4$-web $MW_{a}$ if and only if it is linearizable.
\end{theorem}
\begin{proof}
The necessity is obvious: any $MW_{a}$ is $APW_a$, and by Theorem
10, any $MW_{a}$ is linearizable.

We prove the sufficiency separately for different $a = 1, 2, 3,
4.$

\begin{description}
\item[(i)] \emph{Linearizable $4$-webs $APW_1$}. Note that it
follows from (34), (16), and (18) that for a 4-web $APW_1$ not
only $a_2 = 0$ but also $a_{12} = a_{22} = a_{112} = a_{122} =
a_{222} = 0$. Since by hypothesis our  4-web $APW_1$ is
linearizable, the linearizability conditions (24) and (25) must be
satisfied identically. Substituting
$$
K = K_1 = K_2 = 0, \;\; a_1 = a_{11} = a_{12} = a_{111} = a_{112}
= a_{122} = 0
$$
into (25) and (24), we find that
$$
a_{11} = -\displaystyle \frac{(1-2a)a_1^2 }{a - a^2}
$$
and
$$
a_{111} = -\displaystyle \frac{(6a^2-6a+1)a_1^3 }{(a - a^2)^2}.
$$
Comparing with (48), we see that \emph{a linearizable $4$-web
$APW_1$} is a 4-web $MW_1$. As to the second condition obtained
from (24), it is easy to see that it is a differential consequence
of the first condition.

\item[(ii)] \emph{Linearizable $4$-webs $APW_2$}. Note that it
follows from (34), (16), and (18) that for a 4-web $APW_2$ not
only $a_1 = 0$ but also $a_{11} = a_{12} = a_{111} = a_{112} =
a_{122} = 0$. Since by hypothesis, our  4-web $APW_2$ is
linearizable, the linearizability conditions (24) and (25) must be
satisfied identically. Substituting
$$
K = K_1 = K_2 = 0,\;\; a_1 = a_{11} = a_{12} = a_{111} = a_{112} =
a_{122} = 0
$$
into (24) and (25), we find that
$$
a_{22} = -\displaystyle \frac{a_2^2 }{1 - a}
$$
and
$$
a_{222} = -\displaystyle \frac{a_2^3 }{(1 - a)^2}.
$$
Comparing with (48), we see that \emph{a linearizable $4$-web
$APW_2$} is a 4-web $MW_2$. As to the second condition obtained
from (25), it is easy to see that it is a differential consequence
of the first condition.

\item[(iii)] \emph{Linearizable $4$-webs $APW_3$}. Note that it
follows from (34), (16), and (18) that for a 4-web $APW_3$, the
condition $a_1 = a_2$ implies $a_{11} = a_{12} = a_{22}$ and
$a_{111} = a_{112} = a_{122} = a_{222}$. Since by hypothesis, our
4-web $APW_3$ is linearizable, the linearizability conditions (24)
and (25) must be satisfied identically. Substituting
$$
K = K_1 = K_2 = 0, \;\; a_1 = a_2, a_{11} = a_{12} = a_{22}, \;\;
a_{111} = a_{112} = a_{122} = a_{222}
$$
into (24) and (25), we find that
$$
a_{111} = \displaystyle \frac{4 a_1 a_{11}}{a} - \frac{3
a_1^3}{a^2}
$$
and
$$
a_{111} = -\displaystyle \frac{2 a_1 a_{11}}{a^2} -
\frac{a_1^3}{a^3}.
$$
It follows from these two equations that
$$
a_{11} = -\displaystyle \frac{a_1^2 }{a}.
$$
Comparing this with (48), we see that \emph{a linearizable $4$-web
$APW_3$} is a 4-web $MW_3$. As to the expression for $a_{111}$, it
is easy to see that this expression can be obtained by
differentiation from the expression of $a_{11}$.

\item[(iv)] \emph{Linearizable $4$-webs $APW_4$}. Note that it
follows from (34), (16), and (18) that for a 4-web $APW_4$, the
condition $a_1 = a a_2$ implies
$$
a_{11} = a^2 a_{22} + 2 a_2^2, \;\; a_{12} = a a_{22} + a_2^2
$$
and
$$
\renewcommand{\arraystretch}{1.3}
\begin{array}{ll}
a_{111} = a^3 a_{222} + (7 a^2 a_2 + 2 a a_2) a_{22} + 6 a a_2^3, \\
a_{112} = a^2 a_{222} + 6 a a_2 a_{22} + 2 a_2^3, \\
a_{122} = a a_{222} + 3 a_2 a_{22}.
\end{array}
\renewcommand{\arraystretch}{1}
$$

Since by hypothesis, our  4-web $APW_4$ is linearizable, the
linearizability conditions (24) and (25) must be satisfied
identically. Substituting
$$
K = K_1 = K_2 = 0,\;\; a_1 = aa_2
$$
and the expressions for $a_{11}, a_{12}, a_{111}, a_{112},
a_{122}$ written above into (24) and (25), we find that (25) is
identically satisfied, and (24) reduces to
$$
a_{22} = 0.
$$
As a result, we see from the expressions for $a_{11}$ and $a_{12}$
that they become
$$
a_{11} = 2 a a_2^2, \;\; a_{12} = a_2^2.
$$
Comparing the results with (48), we see that \emph{a linearizable
$4$-web $APW_4$} is a 4-web $MW_4$.
\end{description}
\end{proof}

The next theorem gives another criterion for a 4-web $APW_a$ to be
a 4-web $MW_a$.

\begin{theorem}  A $4$-web $APW_a, \, a = 1, 2, 3, 4,$
is a $4$-web $MW_{a}$ if and only if all its $3$-subwebs are
parallelizable.
\end{theorem}

\begin{proof}
By Theorem 7, for 4-webs $APW_1,\, APW_2$, and $APW_3$, three of
the 3-subwebs are parallelizable. It follows from (38)--(40) that
the fourth 3-subweb of these 4-webs is parallelizable if and only
if we have
$$
  a_{22} = -\displaystyle \frac{a_2^2 }{1 - a}, \;\;
  a_{11} = \displaystyle \frac{(1 - 2 a) a_1^2}{a - a^2},\;\;
a_{11} = a_{12} = a_{22} = \displaystyle\frac{a_1^2}{a},
$$
respectively. Comparing these relations with (48), we see that
  4-webs $APW_1, \linebreak APW_2$, and $APW_3$ are 4-webs
    $MW_1, MW_2$, and $MW_3$, respectively.

For a 4-web $APW_4$, by (41), only the 3-subweb $[1, 2, 3]$ is
parallelizable. But relations (41) show also that other 3-subwebs
are parallelizable if and only if
$$
a_{22} = 0.
$$
But as we already showed in the proof of Theorem 13, the condition
$a_{22} = 0$ implies
$$
a_{11} = 2 a a_2^2, \;\; a_{12} = a_2^2,
$$
and by (48), our 4-web $APW_4$ is a  4-web $MW_4$.
\end{proof}

\textbf{6. 4-webs with a constant basic invariant and equal
curvature forms of their 3-subwebs.} For such 4-webs, we have
\begin{equation}\label{eq:49}
a_1 = a_2 = 0.
\end{equation}
By (14) and (16), conditions (49) imply
\begin{equation}\label{eq:50}
a_{ij} = 0, \;\;\; a_{ijk} = 0, \;\;\;\; i, j, k = 1, 2.
\end{equation}
Now formulas (20)--(22) imply
$$
\crazy{K}{123} = \crazy{K}{124} = \crazy{K}{134} = \crazy{K}{234}
= K
$$
and
$$
\crazy{\Theta}{123} = \crazy{\Theta}{124} = \crazy{\Theta}{134} =
\crazy{\Theta}{234} = K \omega_1 \wedge \omega_2,
$$
i.e., the curvatures and the curvature forms of all 3-subwebs are
equal without vanishing.

We denote such 4-webs by $NW$ (Nakai's webs---see \cite{na96} and
\cite{na98}).

\begin{proposition}  The $4$-webs $NW$  exist, and the general solution
of the system of equations defining such webs depends on one
arbitrary  function of two variables.
\end{proposition}
\begin{proof}
In fact, for 4-webs $NW$, equation (12) becomes $da = 0$. Taking
exterior derivatives of (14), we arrive at the exterior quadratic
equation
$$
\nabla K_1 \wedge \omega_1 + \nabla K_2 \wedge \omega_2 = 0,
$$
where $\nabla K_i = d K_i - 3 K_i \theta$. It is easy to see that
in this case, we have $q = 2$, the Cartan characters $s_1 = 1,\,
s_2 = 1$, and $Q = s_1 + 2 s_2 = 3$.

Next, it follows from the quadratic equation above that
$$
\nabla K_1  = K_{11}  \omega_1 + K_{12} \omega_2, \;\; \nabla K_2
= K_{12}  \omega_1 + K_{22} \omega_2.
$$
Hence  the dimension $N$ of the space of integral elements over a
point equals three,  $N = 3$. Therefore, $Q=N$, and by Cartan's
test, our system of equations is in involution, and its general
solution depends on one arbitrary function of two variables.
\end{proof}

It is easy to see that, in general, a 4-web $NW$ is not
linearizable since conditions (24) and (25) do not hold for such a
web.

However, equations (24) and (25) produce the simple conditions for
linearizability of 4-webs $NW$. Namely, by (49) and (50), it
follows from (24) and (25) that
\begin{equation}\label{eq:51}
K_1 = 0, \;\; K_2 = 0.
\end{equation}

Now we are able to prove the following theorem describing
linearizable 4-webs $NW$.

\begin{theorem}  If a $4$-web $NW$  is linearizable, then
it is parallelizable.
\end{theorem}
\begin{proof}
In fact, for a linearizable 4-web $NW$, we have equations (51).
They and equation (14) imply
$$
dK - 3K \theta = 0.
$$
Taking exterior derivative of this equation and applying (8), we
find that
$$
- 2 K^2 \omega_1 \wedge \omega_2 = 0.
$$
It follows that $K = 0$, and by Theorem 1,  this and (49) prove
that a linearizable 4-web $NW$ is parallelizable.
\end{proof}

Theorem 16 implies the following corollaries.
\begin{corollary}
There exists no linearizable nonparallelizable $4$-webs $NW$.
\end{corollary}

\begin{corollary}
Parallelizable $4$-webs $NW$ exist. They are defined by one
completely integrable equations $da = 0$. The set of such
parallelizable $4$-webs $NW$ depends on one constant $($the
constant basic invariant $a)$.
\end{corollary}
\begin{proof} In fact, since the 4-web $NW$ is parallelizable, we have
$K = K_1 = K_2 \linebreak  = 0$. So, the only remaining equation
defining parallelizable $4$-webs $NW$ is the equation $da = 0, a
\neq 0, 1$. This proves the corollary.
\end{proof}

\vspace*{5mm}

\noindent {\em Author's address}:
Department of Mathematical Sciences, New Jersey Institute  \\
\hspace*{25mm} of Technology,
University Heights,  Newark, N.J. 07102, U.S.A.\\

\noindent
  E-mail address: vlgold@m.njit.edu

\end{document}